\documentclass[psamsfonts]{amsart}
\usepackage{amssymb,amsmath}
\input{diagrams.tex}
\diagramstyle[scriptlabels]
\newcommand{\query}[1]%
{\mbox{}\marginpar{\raggedright\hspace{0pt}{\small\em #1}}}%
\def\ul{\underline}

\def \bbbs{{\mathbb S}}
\def \bbbv{{\mathbb V}}

\def\Hom{\mathop{\rm Hom}\nolimits}
\def\Coker{\mathop{\rm coker}\nolimits}

\def\supp{\mathop{\rm supp}\nolimits}
\def\End{\mathop{\rm End}\nolimits}
\def\id{\mathop{\rm id}\nolimits}

\def\lto{\longrightarrow}

\newcommand{\cE}{{\mathcal E}}
\newcommand{\cF}{{\mathcal F}}
\newcommand{\cH}{{\mathcal H}}
\newcommand{\cK}{{\mathcal K}}
\newcommand{\cG}{{\mathcal G}}

\newcommand{\cM}{{\mathcal M}}
\newcommand{\cO}{{\mathcal O}}

\newcommand{\fm}{{\mathfrak m}}

\def\catp{{\mathfrak P}}
\def\catq{{\mathfrak Q}}
\def\An{{\mathfrak A\mathfrak n}}
\def\Ans{{\mathfrak A\mathfrak n}_S}
\def\Ant{{\mathfrak A\mathfrak n}_T}
\def\sets{{\mathfrak s\mathfrak e\mathfrak t\mathfrak s}}
\def\Equ{\mathop{\rm Equ}\nolimits}
\def\cHom{\mathop{\mathcal Hom}\nolimits}

\theoremstyle{definition}
\newtheorem{dfn}{Definition}[section]

\theoremstyle{plain}
\newtheorem{pro}[dfn]{Proposition}
\newtheorem{thm}[dfn]{Theorem}
\newtheorem{lem}[dfn]{Lemma}
\newtheorem{cor}[dfn]{Corollary}

\theoremstyle{remark}

\begin{document}

\title[Analytic moduli spaces]{Analytic moduli spaces of simple (co)framed
sheaves}

\author{Hubert Flenner}
\address{Fakult\"at f\"ur Mathematik der Ruhr-Universit\"at,
Universit\"atsstr.\ 150, Geb.\ NA 2/72, 44780 Bochum, Germany}
\email{Hubert.Flenner@ruhr-uni-bochum.de}

\author{Martin L\"ubke}
\address{Mathematical Institute, Leiden University, PO box
9512, NL-2300 RA Leiden, The Netherlands}
\email{lubke@math.leidenuniv.nl}

\thanks{This paper was prepared during a visit of the second
author to the University of Bochum which was financed by
EAGER - European Algebraic Geometry Research
Training Network, contract No. HPRN-CT-2000-00099
(BBW 99.0030).}

\begin{abstract}
Let $X$ be a complex space and $\cF$ a coherent $\cO_X$-module. A
{\it $\cF$-(co)framed} sheaf on $X$ is a pair $(\cE,\varphi)$
with a coherent $\cO_X$-module $\cE$ and a morphism of coherent
sheaves
$\varphi : \cF \lto \cE$ (resp. $\varphi : \cE \lto \cF$).
Two such
pairs $(\cE,\varphi)$ and $(\cE^\prime,\varphi^\prime)$ are said
to be {\it isomorphic} if there exists an isomorphism of
sheaves $\alpha : \cE \lto \cE^\prime$ with
$\alpha\circ\varphi = \varphi^\prime$ (resp.
$\varphi^\prime\circ\alpha = \varphi$).
A pair $(\cE,\varphi)$ is called {\it simple} if its only
automorphism is the identity on $\cE$. In
this note we prove a representability theorem
in a relative framework, which implies
in particular that there is a moduli space of simple
$\cF$-(co)framed sheaves on a given compact complex space $X$.
\end{abstract}

\maketitle

\section{Introduction}

Let $X$ be a complex space and $\cF$ a coherent $\cO_X$-module. By a
{\it $\cF$-coframed} sheaf on $X$ we mean a pair $(\cE,\varphi)$
with

(a) $\cE$ is a coherent $\cO_X$-module,

(b) $\varphi : \cF \lto \cE$ is a morphism of coherent
sheaves.

\noindent
(Following \cite{HL1}, \cite{HL2}, a
{\it $\cF$-framed} sheaf is dually a pair $(\cE,\varphi)$ with $\cE$
as above and a morphism $\varphi : \cE \lto \cF$.) Two such
pairs $(\cE,\varphi)$ and $(\cE^\prime,\varphi^\prime)$ are said
to be {\it isomorphic} if there exists an isomorphism of
sheaves $\alpha : \cE \lto \cE^\prime$ with
$\alpha\circ\varphi = \varphi^\prime$.
A pair $(\cE,\varphi)$ is called {\it simple} if its only
automorphism is the identity on $\cE$. The purpose
of this note is to show that there is a moduli space of simple
$\cF$-(co)framed sheaves on a given compact complex space $X$.

More generally, we will show the following relative result. Let
$X\to S$ be a proper morphism of complex spaces. By a family of
$\cF$-coframed sheaves over $S$ (or a $\cF$-coframed sheaf on $X/S$ in
brief) we mean a $\cF$-coframed sheaf $(\cE, \varphi)$ on $X$ that is
$S$-flat. Such a family will be called {\it simple} if its
restriction to each fibre $X(s):=f^{-1}(s)$ is simple.

We consider the set-valued functor
$P : \Ans \lto \sets$ on the category of complex spaces over
$S$ such that
$P(T)$ (for $T \in \Ans$) is the set of all isomorphism
classes of simple $\cF_T$-coframed sheaves $(\cE,\varphi)$
on $X_T/T$, where $X_T:=
X\times_ST$ and $\cF_T$ is the pullback of $\cF$ to $X_T$.
  The main result of this paper is

\begin{thm}\label{1.1}
If $X$ is cohomologically flat over $S$ in dimension 0, then
the functor
$P$ is representable by a (not necessarily
separated) complex space.
\end{thm}

Thus, informally speaking
there is a (relative) moduli space of
$\cF$-coframed sheaves on $X/S$.
An inspection of the proof shows
that the analogous result holds for simple $\cF$-framed sheaves
provided that
$\cF$ is flat over
$S$. The reason that in the case of $\cF$-framed sheaves we need this
additional assumption is that only in this case
the functor $\ul{\Hom}(\cE, \cF)$ is known to be representable (see
\ref{2.3}) whereas  $\ul{\Hom}(\cF, \cE)$ is representable as soon
as $\cE$ is $S$-flat.

Our main motivation for studying moduli spaces of
$\cF$-(co)framed sheaves is the following.
In the case that $\cF$ and $\cE$ are
locally free, a pair
$(\cE,\varphi)$ as above is called a {\it $\cF$-(co)framed vector
bundle} or  {\it holomorphic pair}. Various types of holomorphic
pairs  over compact complex manifolds (e.g. the coframed ones with
$\cF = \cO_X$ or the framed ones with arbitrary $\cF$, see
\cite{OT1}, \cite{OT2}) can be identified with solutions
of so-called {\it vortex equations} via a Kobayashi-Hitchin
type correspondence. On complex surfaces, these
solutions can further be identified with solutions of Seiberg-Witten
equations, and moduli spaces $\cM^{st}$ of {\it stable} holomorphic pairs
(which are open subsets of the moduli spaces $\cM^s$ of simple ones)
can be used to effectively calculate Seiberg-Witten invariants
in several cases.

It is important to notice that the set $\cM^s$ a priori has two analytic
structures. One is given by our result and makes it possible
to determine $\cM^s$ using complex-analytic deformation theory.
The other one,
which is the one relevant in Seiberg-Witten theory,
is given by a gauge theoretical description as in
\cite{LL}. But the main result of that paper is in fact
that these two structures are indeed the same.

Finally we  mention that moduli spaces of stable
$\cF$-framed sheaves on {\it algebraic}
manifolds have been constructed in \cite{HL1}.

Without a (co)framing, a moduli space of simple bundles was
constructed in \cite{KO} and in a more general context in \cite{FS}.
In this paper we follow closely the method of proof in the latter
paper. The main
difficulty is to verify the relative representability of the
general criterion \ref{RCF} for the functor $P$ in \ref{1.1}. For this
we will show in Sect.\ 2 that the functor of endomorphisms of
$\cF$-coframed sheaves is representable. In Sect.\ 3 we show the
openness of the set of points where a coframed sheaf is simple. After
these preparations it will be easy to give in Sect.\ 4 the proof of
\ref{1.1}.

\section{Preparations}

We start with an algebraic lemma.

\begin{lem}
Let $R$ be a ring, $L,M,N$ $R$-modules, and $\alpha : L \lto
M$, $\beta : L \lto N$ $R$-linear maps. If $K := \Coker
((\alpha,-\beta) :L \lto M\oplus N)$, then
$$
\bbbs(K) \cong \bbbs(M)\otimes_{\bbbs(L)}\bbbs(N)
$$
where we consider the symmetric algebras $\bbbs(M)$,
$\bbbs(N)$ as $\bbbs(L)$-algebras via the maps
$\bbbs(\alpha)$, $\bbbs(\beta)$.
\end{lem}

\begin{proof} There are canonical maps $M \lto K$ and $N \lto
K$ given by
$$
m \mapsto \mathrm{residue\ class\ of} (m,0)\ , \ n \mapsto
\mathrm{residue\ class\ of} (0,n).
$$
These maps induce a commutative diagram
\begin{diagram}[s=7mm]
\bbbs(L)&\rTo&\bbbs(M)\\
\dTo&&\dTo\\
\bbbs(N)&\rTo&\bbbs(K)
\end{diagram}
and so, using the universal property of the tensor product,
there is a natural map
$$
\bbbs(M)\otimes_{\bbbs(L)}\bbbs(N) \lto \bbbs(K) .
$$
Conversely, to construct a natural map in the other
direction, note that the first graded piece of
$\bbbs(M)\otimes_{\bbbs(L)}\bbbs(N)$ is just $(M\oplus
N)/L = K$, so there is an induced map
$$
\bbbs(K) \lto \bbbs(M)\otimes_{\bbbs(L)}\bbbs(N) .
$$
It is easy to check that these maps are inverse to each
other. \end{proof}

Now let $S$ be a complex analytic space, and let $\Ans$ be the
category of analytic spaces over $S$. Recall that for a
coherent $\cO_S$-module $\cF$ over $S$ the linear fibre
space $\bbbv(\cF)$ represents the functor
$$
F : \Ans \ni T \mapsto F(T) := \Hom(\cF_T,\cO_T)
$$
(see \cite{Fi} or \cite[II 1.7]{EGA}). Note that $F(T)$ has the
structure of a $\Gamma(T,\cO_T)$-module. Moreover, if $\cG$
is another coherent $\cO_S$-module and $T \mapsto G(T)$ is
the associated functor as above, then a transformation of functors
$F\lto G$ will be called {\it linear} if $F(T) \lto G(T)$ is
$\Gamma(T,\cO_T)$-linear for all $T \in \Ans$. The reader may
easily verify that there is a one-to-one correspondence
between such linear transformations of functors and morphisms
of sheaves $\cG
\lto \cF$.

\begin{pro}
Let $H,F,G : \Ans \lto \sets$ be functors that are
represented by linear fibre spaces
$\bbbv(\cH),\bbbv(\cF),\bbbv(\cG)$, respectively. Let $H
\lto G$ and $F \lto G$ be linear morphisms of functors, and
let $K := H\times_G F$ be the fibered product. Then $K$ is
represented by $\bbbv(\cK)$ with
$$
\cK := \Coker((\alpha,-\beta) : \cG \lto \cH\times\cF),
$$
where $\alpha : \cG \lto \cH$ and $\beta : \cG \lto \cF$
are the morphisms of sheaves corresponding to $H \lto G$ and
$F \lto G$.
\end{pro}

\begin{proof} The spaces $\bbbv(\cH),\bbbv(\cF),\bbbv(\cG)$
are the analytic spectra associated to the symmetric
algebras $\bbbs(\cH),\bbbs(\cF),\bbbv(\cG)$, respectively,
and $H\times_G F$ is represented by
$\bbbv(\cH)\times_{\bbbv(\cG)}\bbbv(\cF)$ which is the
analytic spectrum of
$\bbbs(\cH)\otimes_{\bbbs(\cG)}\bbbs(\cF)$. Hence we need to
verify that there is a natural isomorphism
$$
\bbbs(\cH)\otimes_{\bbbs(\cG)}\bbbs(\cF) \lto \bbbs(\cK),
$$
but this is a consequence of Lemma 2.1. \end{proof}

Let $f : X \lto S$ be a fixed proper morphism of complex
spaces, and let $\cE,\cF$ be coherent $\cO_X$-modules, where
$\cE$ is flat over $S$. Let
$$
H := \underline{\Hom}(\cF,\cE) : \Ans \lto \sets
$$
be the functor given by
$$
H(T) := \Hom_{X_T}(\cF_T,\cE_T),
$$
where $X_T := X\times_S T$ and $\cE_T,\cF_T$ are the
pullbacks of $\cE,\cF$ on $X_T$. We recall the following
fact.

\begin{thm}\label{2.3}
The functor $H$ is representable by a linear fibre space
over $S$.
\end{thm}

For a proof see e.g. \cite[3.2]{Fl2} or \cite{Bi}.

Now let $\varphi : \cF \lto \cE$ and $\varphi^\prime : \cF
\lto \cE^\prime$ be fixed morphisms of coherent sheaves on
$X$. Let us consider the functor
$$
M =
\underline{\Hom}((\cE,\varphi),(\cE^\prime,\varphi^\prime))
: \Ans \lto \sets
$$
defined as follows: For $T \in \Ans$ the elements of $M(T)$
are the pairs
$$
(c,\alpha) \in
\Gamma(\cO_{X_T})\times\Hom_{X_T}(\cE_T,\cE^\prime_T)
$$
such that the diagram
\begin{diagram}[s=7mm,midshaft]
\cF & \rTo^\varphi & \cE\\
\dTo<{c\cdot\id_\cF} & & \dTo>{\alpha}\\
\cF & \rTo^{\varphi^\prime} & \cE^\prime
\end{diagram}
commutes, i.e. such that $c\cdot\varphi^\prime =
\alpha\cdot\varphi$.

\begin{pro}
If $\cO_X$, $\cE$ and $\cE^\prime$ are flat over $S$, then
$M$ is representable by a linear fibre space $\bbbv(\cM)$
over
$S$.
\end{pro}

\begin{proof} By Theorem 2.3 the functors
$$
F := \underline{\Hom}(\cE,\cE^\prime)\ ,\ H :=
\underline{\Hom}(\cO_X,\cO_X)\ ,\ G :=
\underline{\Hom}(\cF,\cE^\prime)\ .
$$
are representable by linear
fibre spaces over $S$. There are natural maps
$$
H \lto G\ ,\ c \mapsto c\cdot\varphi^\prime\ ,
$$
and
$$
F \lto G\ ,\ \alpha \mapsto \alpha\circ\varphi\ .
$$
By definition we have
$$
M = F\times_G H,
$$
so the result follows from Proposition 2.2. \end{proof}

An important property of the sheaf $\cM$ in Proposition 2.4
is given by

\begin{lem}
The following are equivalent.

(a) $\cM$ is locally free.

(b) For every complex space $T\in\Ans$ the canonical map
$$
f_{*}(\cHom((\cE,\varphi),(\cE',\varphi')))\otimes_{\cO_S}
\cO_T
\lto
f_{T*}(\cHom((\cE_T,\varphi_T),(\cE'_T,\varphi'_T)))
$$
is an isomorphism.

Moreover, if one of these conditions holds then
$$
\cM\cong \left[f_{*}(\cHom((\cE,\varphi),(\cE',\varphi')))
\right]^\vee.
\leqno (1)
$$
\end{lem}

\begin{proof}
First note that for every complex space $T\in\Ans$ we have
$$
(\cM_T)^\vee\cong
f_{T*}(\cHom((\cE_T,\varphi_T),(\cE'_T,\varphi'_T))).
\leqno (2)
$$
Applying this to the case $T=S$, (1) follows immediately
from the
assumption that $\cM$ is locally free. Moreover, if (a) is
satisfied we have
$$
\cM_T^\vee\cong\cM^\vee\otimes_{\cO_S}\cO_T\cong
f_{*}(\cHom((\cE,\varphi),(\cE',\varphi')))
\otimes_{\cO_S}\cO_T\,,
$$
where for the last isomorphism we have used (1). Thus (b)
follows.

Conversely, assume that (b) holds.  Using (1) we infer
from the
isomorphism in (b) that
$$
\cHom(\cM_T,\cO_T)\cong \cHom(\cM,\cO_S)\otimes_{\cO_S}\cO_T.
$$
Applying this to $T=\{s\}$, $s\in S$ a reduced point,
it follows that the map
$$
\cHom(\cM,\cO_S)\lto \cHom(\cM,\cO_S/\fm_s)
$$
is surjective for every point $s\in S$. Using standard
arguments (see e.g.\ \cite[7.5.2]{EGA}) we conclude that the functor
$\cHom(\cM,-)$ is exact on the category  of coherent
$\cO_S$-modules whence $\cM$ is locally free, as required.
\end{proof}

\section{Simple $\cF$-coframed sheaves}

As before let $f : X \lto S$ be a proper morphism of complex
spaces, and $\cF$ a fixed $\cO_X$-module. We consider $\cF$-coframed
sheaves $(\cE,\varphi)$ on $X/S$, i.e.\ $\cE$ is a $S$-flat
coherent sheaf on
$X$ and $\varphi : \cF \lto \cE$ is a morphism of
$\cO_X$-modules.

\begin{dfn} $(\cE,\varphi)$ is called {\it simple at} $s\in
S$ if its fibre $(\cE(s), \varphi(s))$ is simple, i.e. if
\begin{diagram}[s=7mm,midshaft]
\cF(s) & \rTo^{\varphi(s)} & \cE(s)\\
\dTo<{c\cdot\id_{\cF(s)}} & &
\dTo>\alpha\\
\cF(s) & \rTo^{\varphi(s)} & \cE(s)
\end{diagram}
is a commutative diagram, then $\alpha = c\cdot\id_{\cE(s)}$.
Moreover, $(\cE,\varphi)$ is said to be {\em simple over $S$} if it
is simple at every point.
\end{dfn}

Notice that this definition of simpleness of $(\cE(s),\varphi(s))$
coincides with the one given in the introduction. Later on we will
will need that the points
$s\in S$ at which $(\cE,\varphi)$ is simple form an open set in
$S$. For this we need the following
considerations.

By Theorem 2.3 and Proposition 2.4 there are coherent
$\cO_S$-modules $\cH$ and $\cG$ such that
$$
\underline{\End}(\cE,\varphi) :=
\underline{\Hom}((\cE,\varphi),(\cE,\varphi))
\quad\mbox{and}\quad
\underline{\End}(\cO_X):=\underline{\Hom}(\cO_X,\cO_X)
$$
are represented by $\bbbv(\cH)$ resp. $\bbbv(\cG)$. Let
$$
\tilde{a} : \cG \lto \cH\quad \mbox{and}\quad
\tilde{b} : \cH \lto \cG
$$
be the $\cO_S$-linear maps associated to the canonical
morphisms of functors
$$
\begin{array}{rclcrcl}
a : \underline{\End}(\cE,\varphi) &\lto&
\underline{\End}(\cO_X)\ ,&\mbox{resp.\ }&
b : \underline{\End}(\cO_X) &\lto&
\underline{\End}(\cE,\varphi)\ ,\\
(c,\alpha)&\longmapsto& c&&
c&\longmapsto&(c,c\cdot\id_\cE)
\end{array}
$$
As $a \circ b = \id_{\underline{\End}(\cO_X)}$ we have
$\tilde{b}\circ\tilde{a} = \id_\cG$. In other words, $\cG$ is
a direct summand of $\cH $ so that $\cH\cong \cG\oplus\cG^\prime$
for some coherent sheaf $\cG'$ on $S$.

\begin{lem}\label{3.2}
The following are equivalent.

(1)  $(\cE,\varphi)$ is simple on $S$.

(2) $\cG^\prime = 0$.

(3)  The canonical
morphism of functors
$b:\underline{\End}(\cO_X) \lto \underline{\End}(\cE,\varphi)$
is an isomorphism.
\end{lem}

\begin{proof}
The functor $\underline{\End}(\cE(s),\varphi(s))$ resp.
$\underline{\End}(\cO_{X(s)})$ on the category $\An$ of all
analytic spaces is represented by
$\cH(s)$ resp. $\cG(s)$. Thus $(\cE(s),\varphi(s))$ is simple if
and only if $\cG(s) \cong \cH(s)$ which is equivalent to the
vanishing of $\cG^\prime(s)$. Using Nakayama's lemma, the
equivalence of (1) and (2) follows.  Finally, the equivalence of (2) and (3) is
immediate from the definition of $\cG'$.
\end{proof}

%
%
%
%

\begin{cor}
If $\cE$ is $S$-flat, then the set of points $s \in S$ at
which $(\cE,\varphi)$ is simple, is an open subset of
$S$.
\end{cor}

\begin{proof}
The set of points $s \in S$ for which
$\cG^\prime/\fm_s\cdot\cG^\prime = 0$ is just the complement
of the support of $\cG^\prime$ and hence Zariski-open in $S$.
Using Lemma \ref{3.2} we get the desired result.
\end{proof}

Recall that a morphism $f : X \lto S$ of complex spaces is
said to be {\em cohomologically flat in dimension} 0 if it is
flat and if for every
$s \in S$ the natural map $f_*(\cO_X) \lto f_*(\cO_{X(s)})$
is surjective.

\begin{cor}
If $(\cE,\varphi)$ is simple and $f : X \lto S$ is
cohomologically flat in dimension 0, then $f_*(\cO_X)$ is a
locally free $\cO_S$-module, and the functor
$\underline{\End}(\cE,\varphi)$ is represented by
$\bbbv(f_*(\cO_X)^\vee)$.
\end{cor}

\begin{proof} The fact that $f_*(\cO_X)$ is locally free
over $S$, is well known (see, e.g.\
\cite[9.7]{FS}). Moreover, since $(\cE,\varphi)$ is simple
we have $\underline{\End}(\cE,\varphi) =
\underline{\End}(\cO_X)$. As
$$
(f_T)_*(\cO_{X_T}) \cong f_*(\cO_X)\otimes_{\cO_S}\cO_T
$$
we get
$$
\Hom_{X_T}(\cO_{X_T},\cO_{X_T}) \cong \Gamma(X_T,\cO_{X_T})
\cong \Hom_T(f_*(\cO_X)^\vee\otimes_{\cO_S}\cO_T,\cO_T)\ ,
$$
so the space $\bbbv(f_*(\cO_X)^\vee)$ represents
$\underline{\End}(\cE,\varphi)$ as desired.
\end{proof}

\section{Proof of Theorem 1.1}

An {\it isomorphism} of two $\cF$-coframed sheaves $(\cE,\varphi)$ and
$(\cE^\prime,\varphi^\prime)$ is an isomorphism
$\alpha : \cE \lto \cE^\prime$ such that
\begin{diagram}[h=4mm,w=7mm,midshaft]
&& \cE \\
&\ruTo^\varphi\\
\cF & & \dTo>\alpha\\
&\rdTo_{\varphi^\prime} \\
&&\cE^\prime
\end{diagram}
commutes. We note that $(\cE,\varphi)$ and
$(\cE^\prime,\varphi^\prime)$ are isomorphic if and only if
there is a pair
$$
(c,\alpha) \in \Gamma(X,\cO_X)\times\Hom_X(\cE,\cE^\prime)
$$
such that

(a) $c$ is a unit in $f_*(\cO_X)$,

(b) $\alpha$ is an isomorphism,

(c) $\alpha\circ\varphi = c\cdot\varphi^\prime$.

\smallskip

Notice that a simple pair $(\cE,\varphi)$ has no automorphism
besides $\id_\cE$. If $S$ is a reduced point then the converse also
holds, i.e.\  $(\cE,\varphi)$ is simple if and only if  $\id_\cE$
is its only automorphism.

\begin{thm}\label{Equ}
Assume that $X$ is cohomologically flat over $S$ in
dimension 0, and let $(\cE,\varphi)$ and
$(\cE^\prime,\varphi^\prime)$ be simple pairs. Then the
functor
$$
F : \Ans \lto \sets\ ,\ F(T) :=
\left\{
\begin{array}{cl}
\{1\} & if\ (\cE,\varphi) \cong
(\cE^\prime,\varphi^\prime),\\
\emptyset & otherwise,
\end{array}
\right.
$$
is representable by a locally closed subspace of $S$.
\end{thm}

\begin{proof}
As the sheaf $f_*(\cO_X)$ is locally free over $S$ we may
assume that it has constant rank, say, $r$ over $\cO_S$. By
Proposition 2.4 the functors
$$
M :=
\underline{\Hom}((\cE,\varphi),(\cE^\prime,\varphi^\prime))\
,\ M^\prime :=
\underline{\Hom}((\cE^\prime,\varphi^\prime),(\cE,\varphi))
$$
are representable by linear fibre spaces $\bbbv(\cM)$ resp.
$\bbbv(\cM^\prime)$, where $\cM$ and $\cM^\prime$ are coherent
$f_*(\cO_X)$-modules. If for some space $T \in \Ans$ the
pairs $(\cE_T,\varphi_T)$ and
$(\cE_T^\prime,\varphi_T^\prime)$ are isomorphic, then by
Corollary 3.4 $\cM_T$ and $\cM_T^\prime$ are locally free
$\cO_T$-modules of rank $r$ on $T$. Thus applying
\cite[9.10]{FS} as in the proof of \cite[9.9]{FS}, we are
reduced to the case that $\cM$ and $\cM^\prime$ are locally
free $\cO_S$-modules of rank $r$. Let us consider the
pairings
$$
\begin{array}{l}
M\times M^\prime \lto \underline{\End}(\cE,\varphi)\ ,\
((c,\alpha),(d,\beta)) \mapsto (cd,\beta\circ\alpha),\\
M^\prime\times M \lto
\underline{\End}(\cE^\prime,\varphi^\prime)\ ,\
((d,\beta),(c,\alpha)) \mapsto (cd,\alpha\circ\beta);
\end{array}$$
these correspond to pairings
\begin{diagram}[s=7mm,midshaft]
\cM^\vee\otimes{\cM^\prime}^\vee & \rTo^\gamma &
f_*(\cO_X),\\ {\cM^\prime}^\vee\otimes\cM^\vee &
\rTo^{\gamma^\prime} & f_*(\cO_X).
\end{diagram}
Using Lemma 2.5 it follows as in the proof of \cite[9.9]{FS}
that our functor $F$ is represented by the open subset
$$
S^\prime :=
S\setminus\supp(\Coker(\gamma)\oplus\Coker(\gamma^\prime)).
$$
\end{proof}

Now we consider the groupoid $\catp \lto \Ans$, where for
$T \in \Ans$ the objects in $\catp(T)$ are the $\cF_T$-coframed
sheaves $(\cE,\varphi)$ on $X_T/T$, where $\varphi : \cF_T \lto \cE$
is $\cO_{X_T}$-linear.  For $(\cE,\varphi) \in \catp(T)$ and
$(\cE^\prime,\varphi^\prime) \in \catp(T^\prime)$, a morphism
$(\cE,\varphi) \lto (\cE^\prime,\varphi^\prime)$ is a pair
$(f,\alpha)$, where $f : T^\prime \lto T$ is an $S$-morphism and
$\alpha : f^*(\cE) \lto \cE^\prime$ is an isomorphism
of coherent  sheaves such that the diagram
\newarrow{Equ} =====
\begin{diagram}[s=7mm]
f^*(\cF_T) & \rTo^{f^*(\varphi)} & f^*(\cE)\\
\dEqu & & \dTo_\alpha\\
\cF_{T^\prime} & \rTo_{\varphi^\prime} & \cE^\prime
\end{diagram}
commutes.

\begin{pro}\label{uvers}${}$

(a) Every object $(\cE_0,\varphi_0)$ in $\catp(s)$, $s\in S$,
admits a
semiuniversal deformation.

(b) Versality is open in $\catp$.
\end{pro}

\begin{proof}
Let $\catq\to\Ans$ be the groupoid where the objects over
a space $T\in \Ans$
are the coherent $\cO_{X_T}$-modules that are $T$-flat.
As usual, given
$\cE\in \catq(T)$ and $\cE'\in \catq(T')$, a morphism
$\cE\to \cE'$ in
$\catq$ consists of a pair $(f,\alpha)$, where $f:T'\to T$
is
an $S$-morphism and $\alpha: f^*(\cE)\to \cE'$ is an
isomorphism of coherent
sheaves. Assigning to a pair $(\cE,\varphi)$ the sheaf
$\cE$ gives a functor
$\catp\to\catq$. It is well known that there are
semiuniversal deformations in
$\catq$ (see \cite{ST} or \cite{BK}) and that versality
is open is $\catq$
(see e.g., \cite{Fl1}).

The fibre of $\catp\to\catq$ over a given object $\cE\in
\catq(T)$ is the groupoid $\catp_\cE\to \Ant$ as explained
in \cite[Sect.\
10]{Bi}. More concretely, given a space $Z\in \Ant$, an
object in $\catp_\cE$
over $Z$ is a morphism
$$
\varphi: \cF\otimes_{\cO_S}\cO_Z\lto \cE\otimes_{\cO_T}\cO_Z.
$$
As the functor underlying $\catp_\cE$ is representable by
Theorem 2.3 we get that
the objects in $\catp_\cE(t)$, $t\in T$, admit semiuniversal
deformations
and that versality is open in $\catp_\cE$.
Applying \cite[10.12]{Bi} gives the
desired conclusion.
\end{proof}

Before proving the main theorem we remind the reader of
the following
criterion for the representability of a functor which we
present for our
purposes in the form as given in \cite[7.5]{FS}; see
also \cite[3.1]{Bi} or
\cite[\S 2]{KO}.

\begin{thm}\label{RCF}
A functor $F:\Ans\to\sets$ is representable by a complex
space over $S$ (resp.
a separated complex space over $S$) if and only if the
following conditions
are satisfied.

{\em (1) (Existence of semiuniversal deformations)}
Every $a_0\in F(s)$, $s\in
S$, admits a semiuniversal deformation.

{\em (2) (Sheaf axiom)} $F$ is of local nature, i.e.\
for every complex space
$T\in\Ans$ the presheaf $T\supseteq U\mapsto F(U)$ on
$T$ is a sheaf.

{\em (3) (Relative representability)} For every
$T\in \Ans$ and $a,\ b\in
F(T)$ the set-valued functor $\Equ(a,b)$ with
$$
\Equ(a,b)(Z):=\left\{
\begin{array}{cl}%
\{ 1\} & if\ a_Z=b_Z,\\
\emptyset & otherwise,
\end{array}%
\right.
$$
is representable by a locally closed (resp.\ closed)
subspace of $T$.

{\em (4) (Openness of versality)} For every $T\in\Ans$
and $a\in F(T)$ the
set of points $t\in T$ at which $a$ is formally versal
is open in $T$.
\end{thm}

\begin{proof}[Proof of Theorem 1.1]
We will verify that the conditions (1)--(4) in
Theorem \ref{RCF} are satisfied. (1) and (4) hold by
Proposition \ref{uvers}. Moreover, (3) is just
Theorem \ref{Equ}. Finally, (2) holds as simple pairs have no
non-trivial automorphism.
\end{proof}


\begin{thebibliography}{FOV}

\bibitem[Bi]{Bi} Bingener, J.: {\em Darstellbarkeitskriterien
f\"ur analytische
Funktoren.} Ann. Sci. \'Ecole Norm. Sup. (4) {\bf 13},
317--347 (1980).

\bibitem[BK]{BK} Bingener, J.: {\em Lokale
Modulr\"aume in der analytischen Geometrie I, II.}
With the cooperation of
S.\ Kosarew. Aspects of Mathematics D2, D3;
Friedr.\ Vieweg \& Sohn,
Braunschweig, 1987.


\bibitem[Fi]{Fi} Fischer, G.: {\em  Complex analytic
geometry.} Lecture Notes
in Math.\ Vol.\ 538. Springer-Verlag, Berlin-New York, 1976.

\bibitem[Fl1]{Fl1} Flenner, H.: {\em Ein Kriterium f\"ur
die Offenheit der
Versalit\"at.}  Math.\ Z.\ {\bf 178}, 449--473 (1981).

\bibitem[Fl2]{Fl2} Flenner, H.: {\em Eine Bemerkung
\"uber relative
Ext-Garben.} Math.\ Ann.\ {\bf 258}, 175--182 (1981).

\bibitem[FS]{FS} Flenner, H.; Sundararaman, D.:
{\em Analytic geometry on
complex superspaces.} Trans.\ AMS {\bf 330}, 1--40 (1992).

\bibitem[EGA]{EGA} Grothendieck, A.; Dieudonn\'e,  J.:
{\em \'El\'ements de  g\'eom\'etrie alg\'ebrique.}
Publ.\ Math.\  IHES {\bf 4,  8,  11,  17,  20,  24,
28,  32},   1961-1967.

\bibitem[HL1]{HL1} Huybrechts, D.; Lehn, M.: {\em Framed
modules and their moduli.} Internat. J. Math. {\bf 6}, 297-324
(1995).

\bibitem[HL2]{HL2} Huybrechts, D.; Lehn, M.: {\em  The
geometry of moduli spaces of sheaves.} Aspects of
Mathematics, E31. Friedr.\ Vieweg \& Sohn, Braunschweig,
1997.

\bibitem[KO]{KO} Kosarew, S.; Okonek, C.: {\em Global
moduli spaces and
simple holomorphic bundles.}
Publ.\ Res.\ Inst.\ Math.\ Sci.\ {\bf  25}, 1--19
(1989).

\bibitem[LL]{LL} L\"ubke, M.; Lupascu, P.:
{\em Isomorphy of the gauge
theoretical and the deformation theoretical
moduli space of simple
holomorphic pairs.} In preparation.

\bibitem[OT1]{OT1}  Okonek, Ch.; Teleman, A.: {\em The
coupled Seiberg-Witten  equations, Vortices, and moduli spaces of
stable pairs.} Internat. J. Math. {\bf 6}, 893--910 (1995).

\bibitem[OT2]{OT2} Okonek, Ch.; Teleman, A.: {\em Gauge theoretical
equivariant Gromov-Witten invariants and the
full Seiberg-Witten invariants of ruled surfaces.}
Preprint (2000).

\bibitem[Ri]{Ri} Rim, D.S.: {\em Formal deformation theory}.
In: S\'eminaire
de G\'eom\'etrie Alg\'ebrique, SGA 7. Lecture Notes in
Mathematics 288,
Springer Verlag Berlin-Heidelberg-New York 1972.

\bibitem[ST]{ST} Siu, Y.T.; Trautmann, G.: {\em
Deformations of coherent
analytic sheaves with compact supports.}
Mem.\ Amer.\ Math.\ Soc.\ {\bf 29} (1981), no.\ 238.

\bibitem[Su]{Su} Suyama, Y.: {\em The analytic moduli
space of simple framed
holomorphic pairs.} Kyushu J.\ Math.\ {\bf  50},
65--82 (1996).


\end{thebibliography}
\end{document}